\documentclass{amsart}
\usepackage{times}
\usepackage{amsfonts}
\usepackage{amsmath}
\usepackage{amscd}
\usepackage{amssymb}
\usepackage{latexsym}
\usepackage{amsthm}

\theoremstyle{plain}

\newtheorem{mainstelling}{Main
Theorem}

\swapnumbers
\newtheorem{theorem}[subsection]{Theorem}
\newtheorem{corollary}[subsection]{Corollary}
\newtheorem{lemma}[subsection]{Lemma}
\newtheorem{proposition}[subsection]{Proposition}

\theoremstyle{definition}
\newtheorem{definition}[subsection]{Definition}

\theoremstyle{remark}

\newtheorem{remark}[subsection]{Remark}

\swapnumbers

\newcommand\NYCCT{\address{Department of Mathematics\\
City University of New York\\
365 Fifth Avenue\\ 
New York, NY 10016 (USA)}
\email{hschoutens@citytech.cuny.edu}}

\usepackage[bookmarks,colorlinks,pagebackref]{hyperref} 

\input{xy}
\xyoption{all}

\newcommand{\emptyprop}{q}

\newcommand \complet[1]{\widehat {#1}}

\newcommand \id{\mathfrak a}

\newcommand \inverse[2]{{#1^{-1}(#2)}}
\newcommand \iso{\cong}
\newcommand \loc{{\mathcal {O}}}
\newcommand \map[1]{{\newcommand{\tmpprop}{#1q}  \if\tmpprop\emptyprop \to\else \xrightarrow{{\phantom{i}{#1}\phantom{i}}}\fi}} 
\newcommand \maxim{\mathfrak m}

\newcommand \pol[2]{#1[#2]}
\newcommand \pow[2]{#1[[#2]]}

\newcommand \pr{\mathfrak p}

\newcommand \range [2]{#1,\dots,#2}
 
\newcommand \restrict [2]{\left.#1\right|_{{#2}}}
\newcommand \rij[2]{(#1_1,\dots,#1_{#2})}

\newcommand \set[2]{\left\{\,#1\mid #2\,\right\}} 
\let\sub\subseteq
\newcommand \sheaf{{\mathcal {F}}}
\newcommand \tensor{\otimes}

\newcommand \op\operatorname

\newcommand{\commdiagram}[9][]{%
\begin{equation}
{\newcommand{\tmpprop}{#1q} 
\if\tmpprop\emptyprop \relax\else \label{#1}\fi}
\begin{aligned}%
\mbox{
\begin{picture}(130,90)%
\put(120,70){\vector( 0,-1){50}}%
\put(10,80){\vector( 1, 0){100}}%
\put(0,70){\vector( 0,-1){50}}%
\put(10,10){\vector( 1, 0){100}}%
\put(115,80){\makebox(0,0)[l]{$#4$}}%
\put(5,80){\makebox(0,0)[r]{$#2$}}%
\put(115,10){\makebox(0,0)[l]{$#9$}}%
\put(5,10){\makebox(0,0)[r]{$#7$}}%
\put(-3,50){\makebox(0,0)[r]{$#5$}}
\put(123,50){\makebox(0,0)[l]{$#6$}}
\put(60,3){\makebox(0,0)[c]{$#8$}}
\put(60,88){\makebox(0,0)[c]{$#3$}}
\end{picture}}
\end{aligned}
\end{equation}}

\newcommand\commtrianglefront[7][]{%
\begin{equation}
{\newcommand{\tmpprop}{#1q} 
\if\tmpprop\emptyprop \relax\else \label{#1}\fi}
\begin{aligned}%
\mbox{
\begin{picture}(120,80)%
\put(55,68){\vector(-1,-2){30}}
\put(65,68){\vector(1,-2){30}}
\put(30,5){\vector(1,0){60}}
\put(60,75){\makebox(0,0)[c]{$#2$}}
\put(25,5){\makebox(0,0)[r]{$#4$}}
\put(95,5){\makebox(0,0)[l]{$#6$}}
\put(60,0){\makebox(0,0)[c]{$#5$}}
\put(37,43){\makebox(0,0)[r]{$#3$}}
\put(83,43){\makebox(0,0)[l]{$#7$}}
\end{picture}}
\end{aligned}
\end{equation}}

\newcommand\commtriangleback[7][]{%
\begin{equation}
{\newcommand{\tmpprop}{#1q}
\if\tmpprop\emptyprop \relax\else \label{#1}\fi}
\begin{aligned}%
\mbox{
\begin{picture}(120,80)%
\put(55,70){\vector(-1,-2){30}}
\put(65,70){\vector(1,-2){30}}
\put(30,5){\vector(1,0){60}}
\put(60,75){\makebox(0,0)[c]{$#2$}}
\put(25,5){\makebox(0,0)[r]{$#6$}}
\put(95,5){\makebox(0,0)[l]{$#4$}}
\put(60,0){\makebox(0,0)[c]{$#5$}}
\put(37,43){\makebox(0,0)[r]{$#7$}}
\put(83,43){\makebox(0,0)[l]{$#3$}}
\end{picture}}
\end{aligned}
\end{equation}}

\newcommand \acf{algebraically closed field}

\newcommand \ch{characteristic}
\newcommand \homo{homomorphism}
\newcommand \CM{Coh\-en-Mac\-au\-lay}
\renewcommand\iff{if and only if}
\newcommand \DVR{discrete valuation ring}

\DeclareMathOperator*{\UP}{ulim}
\newcommand \up[1]{\UP_{#1}}
\newcommand \ul[1]{\seq{#1}\infty}
\newcommand \seq[2]{#1\mathstrut_{#2}}
\newcommand \sr{approximation}

\newcommand  \hull[1]{\mathfrak D(#1)}
\newcommand  \rhull[2]{\mathfrak D_{#1}(#2)}
\newcommand  \ehull[2]{\mathstrut^{#2}\mathfrak D(#1)}
\newcommand  \los{\L os' Theorem}
\newcommand \frob[1]{\mathbf{F}_{#1}}
\newcommand \ulfrob{\frob{}}

\newcommand  \Cech{\v Cech}
\newcommand \ulc{local ultracohomology}
\newcommand \usc{sheaf ultracohomology}

\newcommand \Usc{Sheaf ultracohomology}
\newcommand\nslc[3]{\op{UH}^{#1}_{#2}(#3)}
\newcommand  \nssc[2]{\op{UH}^{#1}(#2, \loc_{#2})}
\newcommand  \BCM[1]{\mathfrak B(#1)}
\newcommand \BS{Brian\c{c}on-Skoda}

\newcommand \gentc[1]{\op{cl}_{\text{gen}}(#1)}

\newcommand  \genrat{generically F-rational}
\newcommand  \genreg{weakly generically F-regular}
\newcommand  \pos[1]{#1^+}
\newcommand \grad[2]{\left [#1\right]_{#2}}
\newcommand \class[3]{ [\frac{#1}{#2}]_{#3}}

\newcommand \zet{\mathbb Z}

\title {Pure subrings of regular rings are pseudo-rational}
\author{Hans Schoutens}
\thanks{Partially supported by a  grant from the National Science Foundation
and a PSC-CUNY grant.}
\NYCCT
\date\today
\subjclass{14B05,13H10,03C20}

\begin{document}

\begin{abstract}   
We prove a generalization of the Hochster-Roberts-Boutot-Kawamata
Theorem conjectured in \cite{SchAsc}: let $R\to S$ be a pure \homo\ of  equi\ch\ zero  Noetherian
local rings. If $S$ is regular, then $R$ is pseudo-rational, and if $R$ is
moreover  $\mathbb Q$-Gorenstein, then it pseudo-log-terminal.
\end{abstract}

\keywords{Tight closure, non-standard Frobenius, rational
singularities, Boutot's Theorem, log-terminal singularities}

\maketitle

\section{Introduction}

Hochster and Roberts showed in \cite{HR}, using finite \ch\ methods, that  
quotient singularities in \ch\ zero are
\CM. This was improved by Boutot in \cite{Bou}
where he shows, using deep vanishing theorems, that they are rational. More
precisely, he shows that if $G$  is the
complexification of a compact Lie group which acts algebraically on an affine
smooth scheme $X$ of finite type over $\mathbb C$, then the quotient $X/G$ has
rational singularities. In algebraic terms, with $X=\op{Spec} B$, this means
that the ring of invariants $A:=B^G$ has rational singularities whenever $B$ is
regular. (In fact, he only requires that $B$ itself has at most rational
singularities,
and there is also a similar result in the analytic category). When $G$ is
finite, Kawamata in \cite{Kaw} showed moreover that $X/G$ has at most
log-terminal singularities,
and the author showed in \cite{SchLogTerm}, using non-standard
tight  closure,  that this remains true for non-finite
$G$, provided $X/G$ is moreover $\mathbb Q$-Gorenstein (a condition that holds
automatically if $G$ is finite).

The goal of the present paper is to extend all these results by removing the
finite type condition. However, since the notion of rational singularities is defined in terms of a
resolution of singularities, which might not be available in such generality,  
we need to replace it by the notion of pseudo-rationality.

\begin{mainstelling}\label{T:psrat}
 Let $A\to B$ be a  cyclically pure \homo\ of
Noetherian   rings containing $\mathbb Q$. If $B$ is regular, then $A$ is
pseudo-rational. 
\end{mainstelling}

 Recall that a \homo\ $A\to  B$ is
\emph{cyclically pure} if $\id=\id B\cap A$ for each ideal $\id$ in $A$;
examples are split, pure or faithfully flat \homo{s}. Since the inclusion $B^G\sub B$ is split
(via the so-called Reynolds operator), Boutot's result is therefore just a
special case of our first main theorem. Theorem~\ref{T:psrat} was
conjectured in
\cite{SchAsc} and proven for algebras of finite
type over an \acf\ in \cite{SchBCM} using canonical big \CM\ algebras. The
analogue in prime \ch\ was proven by Smith in \cite{SmFrat}, but unlike most
applications of tight closure, this proof did not carry over to \ch\ zero, the
reason being that cyclic purity is not preserved under reduction modulo $p$.
In \cite{SchBCM},  we weakened the assumption on $B$ to have Gorenstein  
rational singularities, but at the cost of using a deep theorem due to
Hara: a local $\mathbb C$-algebra $R$ of
finite type has rational singularities \iff\ it is of F-rational type; see \cite{HaRat}.
We showed that this is also equivalent with $R$ being  \genrat.
It is not yet clear whether this remains true in the non-finitely generated
case. To formulate a corresponding generalization in the $\mathbb
Q$-Gorenstein case, we need to make a definition. Call a Noetherian local $\mathbb
Q$-Gorenstein ring $R$ \emph{pseudo-log-terminal}, if its canonical cover $\tilde
R$ (see \S\ref{s:kaw})  is pseudo-rational. In particular, if we are in a category
of local algebras in
which `pseudo-rational' is equivalent with `rational' (e.g., the category of local
algebras essentially of finite type over a field),    
then so is `pseudo-log-terminal' with `log-terminal'  by Theorem~\ref{T:Kaw}.
With this terminology, we get the 
following generalization,  conjectured
in \cite{SchAsc} and proven for algebras of finite type over an \acf\ in
\cite{SchLogTerm}.

\begin{mainstelling}\label{T:logterm} 
Let $R\to S$ be a  cyclically pure
\homo\ of equi\ch\ zero Noetherian local  rings with $S$  regular. If $R$ is   $\mathbb Q$-Gorenstein, then it is pseudo-log-terminal.
\end{mainstelling}

To prove both theorems, we will transform the argument for finitely generated
algebras  given in  \cite{SchLogTerm} by means of the machinery  of faithfully
flat Lefschetz hulls introduced in \cite{SchAsc}. In that paper, we show
that given an equi\ch\ zero Noetherian local ring $R$, we can find a faithfully
flat local $R$-algebra $\hull R$ which is an ultraproduct of rings of prime \ch\
(these latter rings are called \emph{\sr{s}} of $R$ and their ultraproduct is
called a \emph{Lefschetz hull} of $R$). These results enabled us in \cite{SchAsc} to generalize the
alternative 
constructions of tight closure and big \CM\ algebras from the papers
\cite{SchNSTC,SchBCM,SchLogTerm} to arbitrary equi\ch\ zero Noetherian local
rings. Similar applications, although only implicitly using Lefschetz
hulls, can be found in  \cite{SchSymPow,SchBS}.

In the present paper, we will concentrate on one variant coming out of this
work, to wit,  generic tight closure: an element is in the \emph{generic tight
closure} of an ideal if  almost all of its \sr{s} belong to the tight closure of
the corresponding \sr\ of the ideal; see \S\ref{s:gentc} for exact definitions.
Theorem~\ref{T:psrat} will follow from the fact that a \genrat\ ring is
pseudo-rational (see Theorem~\ref{T:genrat}), where we call a ring 
(\emph{generically}) \emph{F-rational} if some ideal generated by a system of parameters is equal to its
(generic) tight closure.   Smith observes in \cite{SmFrat} that F-rationality
in prime \ch\ is equivalent with the top local cohomology of the ring being
Frobenius simple. This enables her to prove that an excellent F-rational
Noetherian
local ring of prime \ch\   is pseudo-regular.  We will not use this result
directly, but rather the method used to prove it. To this end, we also need
Lefschetz hulls for finitely generated algebras over a Noetherian local
ring, as such rings appear in the \Cech\
complex that calculates the local cohomology. This is carried out in
\S\ref{s:RLH}. Therefore, the present proof is entirely self-contained, apart
from  some material taken from \cite{SchAsc}. 

As for Theorem~\ref{T:logterm}, we generalize the notion of an
\emph{ultra-F-regular} local ring introduced in \cite{SchLogTerm} as a
Noetherian local domain $R$ with the property that for each non-zero $c$, we can find an
ultra-Frobenius $\ulfrob^\varepsilon$ such that the morphism $x\mapsto
c\ulfrob^\varepsilon(x)$ is pure (an \emph{ ultra-Frobenius} is an ultraproduct
of powers of Frobenii; see \S\ref{s:uf} below). We then show that the property
of being ultra-F-regular descends under cyclically pure local \homo{s}
(Proposition~\ref{P:cpqreg})  and is preserved under finite extensions which are
\'etale in codimension one (Proposition~\ref{P:et}). Moreover, we show that an
ultra-F-regular local ring is pseudo-rational. 

\subsection*{Open Questions}

\begin{enumerate}
\item Does the converse of Theorem~\ref{T:genrat} also hold, that is to say, is
pseudo-rational equivalent with \genrat? In \cite[Theorem 5.11]{SchBCM}, I gave
a proof of this in the finitely generated case which relies on Hara's result.
\item Does the stronger analogue of Boutot's result also hold, that is to say,
can we weaken the assumption in Theorem~\ref{T:psrat} that $B$ is only
pseudo-rational?  In the finitely generated case,  a proof is available if
$B$ is moreover Gorenstein (\cite[\S5.14]{SchBCM}), but this again depends on
Hara's result.
\item In \cite{SchLogTerm},   using once more Hara's result, it was shown that
for $\mathbb Q$-Gorenstein local domains of finite type over an \acf, the notions
ultra-F-regular and log-terminal are equivalent. Is ultra-F-regular and
pseudo-log-terminal the same for $\mathbb Q$-Gorenstein local domains?
\item Again, we can weaken in the finite type case 
\cite{SchLogTerm} the assumption that $S$ is regular to the assumption that it
is  (pseudo-)log-terminal. Does this also hold in general?
\item For local algebras of finite type over a field of \ch\ zero, rational and
pseudo-rational are the same notions, and so are log-terminal and
pseudo-log-terminal.  For which other categories of equi\ch\ zero Noetherian
local rings is this the case? 
\end{enumerate}

\section{Lefschetz Hulls}\label{s:RLH}

Let $\seq Sw$ be a sequence of rings, where $w$ runs over some infinite set endowed with a
non-principal ultrafilter.  The \emph{ultraproduct} of this sequence is a ring
$\ul S$ given as the homomorphic image of the product $\prod_w\seq Sw$ modulo the
ideal of all sequences which are almost equal to the zero sequence (two sequences $(\seq
aw)$ and $(\seq bw)$ in the product are said to be \emph{almost equal} if $\seq
aw=\seq bw$ for
almost all $w$, that is to say, for all $w$ in some member of the ultrafilter).
When we want to emphasize the index, we denote the ultraproduct $\ul S$ also by
\begin{equation*}
\up w\seq Sw
\end{equation*}
and similarly, the image of a sequence $(\seq aw)$ in $\ul S$ is denoted  $\up
w\seq aw$ or simply $\ul a$.
In case all rings are equal, say $\seq Sw:=S$, their
ultraproduct is called an  \emph{ultrapower} of $S$.  For more details, see   \cite[\S9.5]{Hod} or
\cite{EkUP}, or the brief review in \cite[\S2]{SchNSTC}.

\subsection{Lefschetz hulls} Let $K$ be an uncountable \acf\ of \ch\
zero. In \cite{SchAsc}, we associate to every  Noetherian local ring 
$R$ whose residue field is contained in $K$, a  faithfully flat
Lefschetz hull, that is to say, a faithfully flat extension $R\sub  \hull
R$ such that $\hull R$ is an ultraproduct of prime \ch\ (complete) Noetherian
local rings $\seq Rw$. Any sequence of prime \ch\ complete Noetherian local rings
$\seq Rw$
whose ultraproduct is equal to $\hull R$ is called an \emph{\sr} of $R$.  
For the extent to which the assignment $R\mapsto \hull R$ is functorial, we
refer to the cited paper. All we need in the present paper is that if
$R\to S$ is a local \homo\ of Noetherian local rings whose residue field is
contained in $K$, then there is a \homo\ $\hull R\to\hull S$ making the
following diagram commute 
\commdiagram R {} {\hull R} {} {} S {} {\hull S.}

For the remainder of this section, $R$ is an equi\ch\ zero
Noetherian local ring, $\seq Rw$ an \sr\ of $R$ and $\hull R$ its
Lefschetz hull.
For each $w$, let $\frob w$ denote the
Frobenius on $\seq Rw$, that is to say the \homo\ given by $x\mapsto
x^{p(w)}$, where $p(w)$ is the \ch\ of $\seq Rw$.
Given a positive integer $e$, let $\mathstrut^e
\seq Rw$ denote
 the $\seq Rw$-algebra
structure on $\seq Rw$ given by $\frob w^e$.   It follows that $\frob w^e\colon
\seq Rw\to
\mathstrut^e \seq Rw$ is $\seq Rw$-linear.

\subsection{Ultra-Frobenius}\label{s:uf} 
A \emph{non-standard integer} is an element
$\varepsilon$ of the ultrapower $\ul\zet$ of $\zet$, that is to say, an
ultraproduct of integers $\seq ew$. If almost all $\seq ew$ are positive, then we call
$\varepsilon$ \emph{positive}. For each positive non-standard integer $\varepsilon$, let 
$\ulfrob^\varepsilon\colon R\to \hull R$ be the ultraproduct of the $\frob
w^{\seq ew}$, that is to say, for $x\in R$ with \sr\ $\seq xw$, we have
\begin{equation*}
\ulfrob^\varepsilon(x):= \up w \frob w^{\seq ew}(\seq xw)\in \hull R.
\end{equation*} 
As in \cite{SchLogTerm}, we will call any \homo\ $R\to \hull R$ of the form
$\ulfrob^\varepsilon$ for some $\varepsilon$ an \emph{ultra-Frobenius}. If
$\varepsilon=1$, then the corresponding ultra-Frobenius is just the
\emph{non-standard Frobenius} introduced  in \cite{SchAsc}.

For each positive non-standard integer $\varepsilon$, we may view $\hull R$ as
an $R$-algebra via the \homo\ $\ulfrob^\varepsilon$. To denote this algebra
structure, we will write $\ehull R\varepsilon$ (in   \cite{SchLogTerm}, the
alternative notation $(\ulfrob^\varepsilon)_*\hull R$ was used). In other words,
the $R$-algebra structure on $\ehull R\varepsilon$ is given by $x\cdot \alpha:=
\ulfrob^\varepsilon(x)\alpha$, for $x\in R$ and $\alpha\in\hull R$.

One of the major drawbacks of the functor $\mathfrak D$ is its local nature. In
particular, since a localization $R\to R_\pr$ is not a local \homo, there is no
obvious map from $\hull R$ to $\hull{R_\pr}$. Below we will have to deal with
localizations of the form $R_y$, and hence we need a notion of Lefschetz
hull for such (non-local) rings as well.

\subsection{$R$-\sr{s}}  Let $Y$ be a tuple of
indeterminates and let $f\in \pol RY$, say of the form $f=\sum_{\nu\in N} a_\nu
Y^\nu$ with $a_\nu\in R$ and $N$ a finite index set. If $\seq{a_\nu}w$ is an
\sr\ of $a_\nu$, for each $\nu\in N$, then we call the sequence of polynomials
 $\seq fw:=\sum_{\nu\in N}\seq{a_\nu}wY^\nu$ an \emph{$R$-\sr} of $f$.

One checks that any two $R$-\sr{s} of a polynomial $f$  are almost equal.
Similarly, if $I:=\rij fs$ is an ideal in $\pol RY$ and $\seq{f_i}w$ is an
$R$-\sr\ of $f_i$, for each $i$, then  we call the sequence $\seq
Iw:=(\seq{f_1}w,\dots,\seq{f_s}w)\pol{\seq Rw}Y$  an \emph{$R$-\sr} of $I$, and
if $S=\pol RY/I$, then we call the sequence $\seq Sw:=\pol{\seq Rw}Y/\seq Iw$ an
\emph{$R$-\sr} of $S$. 

\subsection{Relative hulls} If $S$ is a finitely generated $R$-algebra and $\seq
Sw$ is an $R$-\sr\ of $S$, then the ultraproduct of the $\seq Sw$ is called the
\emph{(relative) $R$-hull} of $S$ and is denoted $\rhull RS$.

If $\pol RZ/J$ is another presentation of $S$ as an $R$-algebra, then we have
substitution maps $Y\mapsto \mathbf a$ and $Z\mapsto \mathbf b$ which induce
isomorphisms modulo $I$ and $J$ respectively, where $\mathbf a$ and $\mathbf b$
are tuples of polynomials in the $Z$ and $Y$-variables respectively. Let $\seq
{\mathbf a}w$ and $\seq {\mathbf b}w$  be $R$-\sr{s} of these respective tuples
and let $\seq Jw$ be an $R$-\sr\ of $J$. By \los\   the substitutions $Y\mapsto
\seq {\mathbf a}w$ and $Z\mapsto \seq {\mathbf b}w$ induce for almost all $w$
isomorphisms modulo $\seq Iw$ and $\seq Jw$ respectively. It follows that the
ultraproduct of the $\pol{\seq Rw}Y/\seq Iw$ is isomorphic to the ultraproduct
of the $\pol{\seq Rw}Z/\seq Jw$, showing that $\rhull RS$ is independent from
the particular presentation of $S$ and from the particular choice of $R$-\sr{s}.

Since $\rhull RS$ is naturally a $ \hull R$-algebra and since by \los\ the tuple
$Y$ is algebraically independent over $\hull R$, we get a natural $\pol{\hull
R}Y$-algebra structure, whence an $\pol RY$-algebra structure, on $\rhull RS$.
Under the natural \homo\ $\pol RY\to \rhull RS$, we get $I\rhull RS=0$,
so that this induces a
\homo\
$S\to
\rhull RS$,
endowing $\rhull RS$ with a canonical $S$-algebra structure. We can now extend
the notion of $R$-\sr\ of an element $a$ or an ideal $\id$ in a finitely
generated $R$-algebra $S$ as follows. Let $S:=\pol RY/I$ and choose a polynomial
$f\in\pol RY$ and an ideal $\mathfrak A$ in $\pol RY$ so that their
images in $S$ are  respectively $a$ and $\id$. Let $\seq fw$, $\seq{\mathfrak
A}w$ and $\seq Sw$ be $R$-\sr{s} of $f$, $\mathfrak A$ and $S$ respectively. We call the image
$\seq aw$ of $\seq fw$ in $\seq Sw$ (respectively, the ideal $\seq\id w:=\seq{\mathfrak
A}w\seq Sw$) an \emph{$R$-\sr} of $a$ (respectively, of $\id$). Note that the
ultraproduct of the $\seq aw$ (respectively, of the $\seq\id w$) is equal to the
image of $a$ in $\rhull RS$ (respectively, equal to the ideal $\id\rhull RS$),
showing that any two $R$-\sr{s} are almost equal.

If $S\to T$ is an $R$-algebra \homo\ of finite type, then this extends to an
$R$-algebra \homo\ $\rhull RS\to\rhull RT$ giving rise to a commutative diagram
\commdiagram [rhull] S {} {\rhull RS} {} {} T {} {\rhull RT.} In particular,
$\rhull R\cdot$ is a functor on the category of finitely generated
$R$-algebras. 
The argument is the same as in \cite[\S3.2.4]{SchNSTC} and we leave the details
to the reader.

\section{Generic Tight Closure}\label{s:gentc}

One of the tight closure notions introduced in \cite{SchAsc} is generic tight
closure. In this section, we review the definition and (re)prove some of its
main properties. As above, fix an uncountable \acf\ $K$ of \ch\ zero. Throughout
this section,  $(R,\maxim)$  will denote  a Noetherian local ring with residue
field contained in $K$ and  $(\seq Rw,\seq\maxim w)$ one of its \sr{s}.  For
generalities on (\ch\ $p$) tight closure, see \cite{HuTC}.

\begin{definition}
An element $z\in R$ lies in the \emph{generic
tight closure} of an ideal $\id\sub  R$, if almost all $\seq zw$ lie in the tight closure
$\seq\id w^*$ of $\seq \id w$, where $\seq zw$ and $\seq \id w$ are \sr{s} of
$z$ and $\id$ respectively.  
\end{definition}

We denote the generic tight closure of an ideal
$\id$ by $\gentc\id$. One easily checks that
\begin{equation}
\label{eq:gentc}
\gentc\id = (\up w \seq\id w^*)\cap R
\end{equation}
where the contraction of course is with respect to the canonical embedding $R\to
\hull R$.
It follows that $\gentc\id$ is an ideal, containing $\id$, with the property
that $\gentc{\gentc\id}=\gentc\id$. We say that an ideal $\id$ is
\emph{generically tightly closed} if  $\id=\gentc\id$. The proof of the following easy fact is left to the reader. 

\begin{lemma}\label{L:col} If $\id\sub  R$ is a generically tightly closed
ideal, then so is any colon ideal $(\id:_R\mathfrak b)$, for $\mathfrak b\sub 
R$. \qed
\end{lemma}

\begin{theorem}\label{T:reg} If $R$ is regular, then every ideal is generically
tightly closed.
\end{theorem}
\begin{proof} By \cite[Theorem 5.2]{SchAsc}, almost all $\seq Rw$ are regular,
and hence all ideals in $\seq Rw$ are tightly closed by \cite[Theorem
1.3]{HuTC}.
The assertion then follows from \eqref{eq:gentc} and faithful flatness.
\end{proof}

\begin{theorem}[Persistence]\label{T:per} If $R\to S$ is a local \homo\ and
$\id$ an ideal in $R$, then $\gentc\id S\sub  \gentc{\id S}$.
\end{theorem}
\begin{proof} Immediate from \eqref{eq:gentc} and the fact  that persistence
holds for each $\seq Rw\to \seq Sw$, where $\seq Sw$ is an \sr\ of $S$ (note
that $\seq Rw$ is complete, so that \cite[Theorem 2.3]{HuTC} applies).
\end{proof}

\begin{theorem}[Strong Colon Capturing]\label{T:CC} Let  $\rij xd$ be part of a
system of parameters of $R$. For each $i$, the element $x_i$ is a non-zero
divisor modulo $\gentc{\rij x{i-1}R}$.
\end{theorem}
\begin{proof} By downward induction on $i$, it suffices to prove the assertion
for $i=d$. To this end, suppose $zx_d\in\gentc I$ with $I:=\rij x{d-1}R$. Let
$\seq Rw$, $\seq zw$ and $\seq{x_i}w$ be   \sr{s} of $R$, $z$ and $x_i$
respectively and put $\seq Iw:=(\seq{x_1}w,\dots,\seq{x_{d-1}}w)\seq Rw$. By
\cite[Corollary 5.3]{SchAsc}, almost all $(\seq{x_1}w,\dots,\seq{x_d}w)$ are
part of a system of parameters in $\seq Rw$ and  $\seq zw\seq{x_d}w\in \seq
Iw^*$. Since  each $\seq Rw$ is complete, Strong Colon Capturing holds for it,
that is to say, $\seq{x_d}w$ is a non-zero divisor modulo $\seq Iw^*$ (see 
\cite[Theorem 3.1A and Lemma 4.1]{HuTC}). Therefore, $\seq zw \in \seq Iw^*$,
whence $z\in\gentc I$, as we needed to show.
\end{proof}

\begin{remark} In particular, the usual Colon Capturing holds, that is to say,
for each $i$, we have an inclusion $(\rij x{i-1}R:x_i)\sub  \gentc{\rij
x{i-1}R}$. The same proof can also be used to prove the following stronger
version (compare with \cite[Theorem 9.2]{HuTC}):
let $\pol \zet X\to R$ be given by $X_i\mapsto x_i$ and let $I,J\sub \pol\zet
X$ be ideals. We have an inclusion
\begin{equation}\label{eq:strcc}
(\gentc{IR}:JR) \sub \gentc{(I:J)R}.
\end{equation}
\end{remark}

\begin{corollary}\label{C:regseq} 
If $\rij xd$ is part of a system of parameters
in $R$ and if $\rij xdR$ is generically tightly closed, then so is each   $\rij
xiR$, for $i=\range 1d$. In particular, $\rij xd$ is a regular sequence.
\end{corollary}
\begin{proof} The last assertion is clear from Colon Capturing and the first
assertion. For the first assertion, it suffices to treat the case $i=d-1$, by
downwards induction on $i$. Let $I:=\rij x{d-1}R$ and let $z\in\gentc I$.
Clearly, $z\in\gentc {I+x_dR}$ and this latter ideal is just $I+x_dR$ by
hypothesis. Write $z=a+rx_d$, with $a\in I$ and $r\in R$. Therefore,
$z-a=rx_d\in\gentc I$. Since $x_d$ is a non-zero divisor modulo $\gentc I$ by 
Theorem~\ref{T:CC}, we get $r\in\gentc I$. So, we proved that $\gentc
I=I+x_d\gentc I$.   Nakayama's Lemma then yields that $I=\gentc I$.
\end{proof}

\begin{theorem}[\BS]\label{T:BS} The generic tight closure of  an ideal
$\id\sub  R$ is contained in its integral closure.  If $\id$ is generated by
at most $n$ elements, then the integral closure of $\id^{m+n}$ is contained in
$\gentc{\id^{m+1}}$, for each $m$.
\end{theorem}
\begin{proof} Let $z\in\gentc\id$. In  order to show that $z$ is integral over
$\id$, it suffices by \cite[Lemma 2.3]{HoTIC} to show that $z\in\id V$, for each
\DVR\ $V$ such that $R\to V$ is a local \homo. Now, persistence (Theorem~\ref{T:per}) yields that $z$ lies in $\gentc{\id
V}$, whence, by Theorem~\ref{T:reg},  in $\id V$.

Assume next that $z$ lies in the integral closure of $\id^{m+n}$, for some $m$
and for $n$ the number of generators of $\id$. Taking an integral equation
witnessing this fact and considering \sr{s}, we see that almost all $\seq zw$
lie in the integral closure of $\seq\id w^{m+n}$, where $\seq zw$ and $\seq\id
w$ are \sr{s} of $z$ and $\id$ respectively.  By the tight closure \BS\ Theorem
(see for instance \cite[Theorem 5.7]{HuTC} for an easy proof), 
almost all $\seq zw$ lie in the tight closure of $\seq \id w^{m+1}$ and the
result follows.
\end{proof}

\subsection{Comparison with other tight closure operations}\label{s:tc} By
\cite[Theorem 6.21]{SchAsc}, the generic tight closure of an ideal $\id$ is
contained in its non-standard tight closure, provided $R$ is analytically
unramified. This latter condition is imposed to insure the existence of uniform
test elements (\cite[Proposition 6.20]{SchAsc}).

If $R$ is moreover equidimensional and universally catenary, then by
\cite[Proposition 7.13]{SchAsc}, the $\mathfrak B$-closure $\id\BCM R\cap R$ of
$\id$ is contained in its generic tight closure, with equality if $\id$ is
generated by a system of  parameters. Here $\BCM R$ denotes the canonical big
\CM\ algebra associated to $R$ from \cite[\S7]{SchAsc}. (In the special case
that $R$ is a complete domain with algebraically closed residue field, $\BCM R$
is obtained as the ultraproduct of the absolute integral closures $\seq Rw^+$.)

\section{Generic  F-Rationality}

As before, $R$ is a Noetherian local ring with residue field contained in the
fixed uncountable \acf\ $K$ and $\seq Rw$ is an \sr\ of $R$.

\begin{definition} We say that $R$ is \emph{\genrat}, if there exists a system
of parameters $\mathbf x$ in $R$ such that $\mathbf xR$ is generically tightly
closed.
\end{definition}

Let us say that $R$ is \emph{$\mathfrak B$-rational}, if there exists a system of
parameters $\mathbf x$ such that $\mathbf xR=\mathbf x\BCM R\cap R$. We will
prove below that a ring is \genrat\ \iff\ its completion $\complet R$ is. We
leave it as an exercise to prove that the same property with `$\mathfrak
B$-rational' instead of `\genrat' also holds. Therefore, in view of our
discussion in \S\ref{s:tc}, a ring is \genrat\ \iff\ it is $\mathfrak
B$-rational.

\begin{theorem}\label{T:CM} If $R$ is \genrat, then it is  \CM.
\end{theorem}
\begin{proof} Let $\mathbf x$ be a system of parameters in $R$ such that
$\mathbf xR$ is generically tightly closed.  By Corollary~\ref{C:regseq}, the
sequence $\mathbf x$  is regular and hence $R$ is \CM. 
\end{proof}

\begin{theorem}\label{T:norm} 
If $R$ is \genrat, then any ideal  generated by
part of a system of parameters is generically tightly closed. In particular, $R$
is normal.
\end{theorem}
\begin{proof} By Theorem~\ref{T:CM}, we know that $R$ is \CM. By
Corollary~\ref{C:regseq}, it suffices to show that any ideal generated by a
system of parameters $\rij yd$ is generically tightly closed.  Reasoning on the
top local cohomology,  we can find $t\geq 1$ and $a\in R$ such that  $\rij
ydR=((x_1^t,\dots,x_d^t)R:_R a)$ (see for instance the proof of \cite[Lemma
4.1]{HuTC}). Therefore, if we can show that $(x_1^t,\dots,x_d^t)R$ is 
generically tightly closed, then so will  $\rij ydR$ be by Lemma~\ref{L:col}.
Hence we have reduced to the case that $y_i=x_i^t$, for some $t\geq 1$. 

Let $z\in\gentc{(x_1^t,\dots,x_d^t)R}$. We need to show that
$z\in(x_1^t,\dots,x_d^t)R$. If some $zx_i$ does not lie in
$(x_1^t,\dots,x_d^t)R$, we may replace our original $z$ by this new element.
Therefore, we may assume that 
\begin{equation*} 
z\rij xdR\sub  (x_1^t,\dots,x_d^t)R.
\end{equation*} 
Since $\rij xd$ is $R$-regular, we have 
\begin{equation*} 
((x_1^t,\dots,x_d^t)R:\rij xdR)=(x_1^t,\dots,x_d^t, y^{t-1})R,
\end{equation*} 
where $y:=x_1\cdots x_d$. In summary, we may assume that
$z=uy^{t-1}$, for some $u\in R$. By \eqref{eq:strcc}, we then get
\begin{align*}
u\in (\gentc{(x_1^t,\dots,x_d^t)R}:y^{t-1}) &\sub
\gentc{((x_1^t,\dots,x_d^t)R:y^{t-1})}\\
&=\gentc {\rij xdR}=\rij xdR.
\end{align*}
Therefore, $z=uy^{t-1}$ lies in $(x_1^t,\dots,x_d^t)R$, as we wanted
to show.
 
In order to prove that $R$ is normal, if suffices to show that any height one
principal ideal $aR$ is integrally closed. Since the integral closure of $aR$ is
contained in $\gentc{aR}$ by Theorem~\ref{T:BS}, and since $a$ is part of a
system of parameters, the conclusion follows from the first assertion.
\end{proof}

\begin{proposition}\label{P:comp} A local ring   $R$ is \genrat\ \iff\  its
completion $\complet R$ is. In particular, a \genrat\ ring is analytically
unramified.
\end{proposition}
\begin{proof} Let $\mathbf x$ be a system of parameters in $R$ such that
$\mathfrak n:=\mathbf xR$ is generically tightly closed.  I claim that
$\mathfrak n\complet R$ is generically tightly closed, from which it follows
that $\complet R$ is \genrat. To this end, let $\complet z\in \complet R$ be in
the generic tight closure of $\mathfrak n\complet R$.   Write $\complet z=z+
\complet w$ with $z\in R$ and $\complet w\in \mathfrak n\complet R$. It follows
that $z\in \gentc{\mathfrak n\complet R}$. Let $J$ be the ultraproduct of the
$\seq {\mathfrak n}w^*$, where $\seq {\mathfrak n}w$ is an \sr\ of $\mathfrak
n$.  Since $\seq Rw$ is also an \sr\ for $\complet R$, we get $\gentc{\mathfrak n\complet R}=J\cap \complet
R$ by \eqref{eq:gentc}. Hence
$z\in J$, and since $J\cap R=\gentc{\mathfrak n}=\mathfrak n$, we get
$\complet z=z+\complet w\in\mathfrak n\complet R$. 

Conversely, suppose $\complet R$ is \genrat. Let $\mathbf x$ be a system of
parameters in $R$. Since $\mathbf x$ is then also a system of parameters in
$\complet R$, the ideal $\mathbf x\complet R$ is generically tightly closed by
Theorem~\ref{T:norm}. Let $a\in\gentc{\mathbf xR}$. By  Theorem~\ref{T:per}, we
get that $a$ lies in the generic tight closure of $\mathbf x\complet R$ whence
in $\mathbf x\complet R$. By faithful flatness, we get $a\in\mathbf xR$, proving
that $R$ is \genrat.

To prove the last assertion, assume $R$ is \genrat. Hence so is $\complet R$ by
what we just proved. Therefore, $\complet R$ is normal by Theorem~\ref{T:norm},
whence a domain, showing  that $R$ is analytically unramified.
\end{proof}

\begin{corollary}\label{C:srCM} If $R$ is \genrat, then almost all $\seq Rw$ are
\CM\ and normal.
\end{corollary}
\begin{proof} Since $\complet R$ and $R$ have the same \sr{s}, we may assume by
Proposition~\ref{P:comp} that $R$ is complete.   Theorems~\ref{T:CM} and
\ref{T:norm} yield that $R$ is normal and \CM. By \cite[Theorem 5.2]{SchAsc},
almost all $\seq Rw$ are \CM. Since $R$ satisfies Serre's condition $(\op R_1)$,
so do almost all $\seq Rw$ by \cite[Theorem 5.6]{SchAsc}. Together with the fact
that almost all $\seq Rw$ are \CM, we get from Serre's criterion for normality
(see for instance \cite[Theorem 23.8]{Mats}) that almost all $\seq Rw$ are
normal.
\end{proof}

\begin{proposition}\label{P:Frat} 
If almost all $\seq Rw$
are F-rational, then $R$ is \genrat. The converse holds if $R$ is moreover
Gorenstein.
\end{proposition}
\begin{proof} Let $\mathbf x$ be a system of parameters in $R$, with \sr\ $\seq
{\mathbf x}w$, and let $z$ be in the generic tight closure of $\mathbf xR$. By
\cite[Corollary 5.4]{SchAsc}, almost all $\seq{\mathbf x}w$ are systems of
parameters in $\seq Rw$. Hence, by definition of F-rationality, $\seq{\mathbf
x}w\seq Rw$ is tightly closed. Therefore, if $\seq zw$ is an \sr\ of $z$, then
$\seq zw\in\seq{\mathbf x}w\seq Rw$. Taking ultraproducts, we see that $z$ lies
in $\mathbf x\hull R$ and hence by faithful flatness, in $\mathbf xR$, showing
that $R$ is \genrat.

Suppose next that $R$ is Gorenstein and \genrat. Towards a contradiction, assume
almost all $\seq Rw$ are not F-rational. If $J$ is the ultraproduct of the
$(\seq{\mathbf x}w\seq Rw)^*$, then this means that $\mathbf x\hull
R\varsubsetneq J$. On the other hand, by \eqref{eq:gentc} and our assumption,
$J\cap R= \mathbf xR$. Put $S:=R/\mathbf xR$. By \cite[\S4.9]{SchAsc}, we have
an isomorphism $\hull S\iso\hull R/\mathbf x\hull R$ and $\hull S$ is an
ultrapower of $S\tensor_kK$, where $k$ is the residue field of $R$. Since $S$ is
Gorenstein, so is $S\tensor_kK$, whence also $\hull S$, since   the Gorenstein
property is
first order definable (see for instance \cite{SchEC}). Let $a\in R$ be such that its image in $S$ generates the socle
of this ring. By faithful flatness, $a$ is a non-zero element in the socle of
$\hull S$, whence must generate it. Since $JS\neq 0$, we must have $a\in J$
whence $a\in J\cap R=\mathbf x R$, contradiction.
\end{proof}

\begin{remark}
 Note that by Smith's result \cite[Theorem 3.1]{SmFrat}, an
F-rational excellent local ring   is pseudo-rational; the converse holds by
\cite{HaRat}. It follows that if almost all \sr{s} of $R$ are pseudo-rational,
then $R$ is \genrat, whence pseudo-rational by Theorem~\ref{T:genrat} below.  I
do not know whether the converse also holds.
\end{remark}

Let us call $R$ \emph{\genreg}, if each ideal $\id\sub  R$ is generically
tightly closed.  By Theorem~\ref{T:reg}, any regular local ring is \genreg.
By a similar argument as in the proof of Proposition~\ref{P:comp}, one can
show that $R$ is \genreg\ \iff\ its completion is. If
a ring  is \genreg, then it is \genrat; the converse is true for Gorenstein rings,
as we now prove.

\begin{theorem}
 If $R$ is Gorenstein and \genrat, then it is \genreg.
\end{theorem}
\begin{proof} Given an arbitrary  ideal $\id\sub  R$, we need to show that
$\id=\gentc\id$. Since $\id$ is the intersection of $\maxim$-primary ideals, we
easily reduce to the case that $\id$ is $\maxim$-primary. Choose a system of
parameters $\mathbf x$ such that $\mathbf xR\sub  \id$. By
Theorem~\ref{T:norm}, the ideal $\mathbf xR$ is generically tightly closed.
Since $R$ is Gorenstein,  
\begin{equation*}
\id =(\mathbf xR:(\mathbf xR:\id))
\end{equation*} which is a generically tightly closed ideal by
Lemma~\ref{L:col}.
\end{proof}

\begin{proposition}\label{P:cp} 
Let $R\to S$ be a cyclically pure, local \homo\
between   Noetherian local rings whose residue fields lie in $K$. If $S$ is
\genreg, then so is $R$.
\end{proposition}
\begin{proof} Let $z\in\gentc\id$, for $\id$ an ideal in $R$. By
Theorem~\ref{T:per}, the image of $z$ in $S$ lies in the generic tight closure
of $\id S$, which by assumption is just $\id S$. Hence $z\in \id S\cap R=\id$.
\end{proof}

\begin{remark} It is well-known that the localization of an F-rational ring is
again F-rational (see \cite[Theorem 4.2]{HuTC}; the same property for weakly
F-regular rings though is still open). However, since   Lefschetz hulls are
not compatible with localization, I do not know whether the localization of a
\genrat\ ring is again \genrat.
\end{remark}

The next \BS\ type theorem was proven first in \cite{LT} for pseudo-rational
local rings. Since we will show in the next section that a \genrat\ local ring
is pseudo-rational, this version generalizes their result.

\begin{theorem} 
If $R$ is a $d$-dimensional \genrat\ local ring, then the
integral closure of $\id^{m+d}$ is contained in $\id^{m+1}$, for all $m$ and all
ideals $\id\sub  R$.
\end{theorem}
\begin{proof} We follow the argument in \cite[Theorem 6.4]{SchBCM}, where the
special case that $R$ is of finite type over an \acf\ is proven. Let $a$ be an
element of the integral closure of $\id^{m+d}$. Assume first that $\id$ is
generated by a system of parameters. Therefore,  $a$ lies in
$\gentc{\id^{m+1}}$, by Theorem~\ref{T:BS},  whence in $\id^{m+1}$, by
Lemma~\ref{L:pow} below. This proves the assertion for parameter ideals. Assume
next that $\id$ is merely $\maxim$-primary, where $\maxim$ is the maximal ideal
of $R$. In that case, $\id$ admits a reduction $I$ generated by a system of
parameters. Since  $I^{m+d}$ is then a reduction of  $\id^{m+d}$, we get that 
$a$ lies in the integral closure of $I^{m+d}$, whence in $I^{m+1}$, by the first
case, and, therefore, ultimately in  $\id^{m+1}$, also establishing this case.
For arbitrary $\id$,
write $\id$ as the intersection of all $\id+\maxim^n$ and use the previous case.
\end{proof}

\begin{lemma}\label{L:pow} 
If $R$ is a \genrat\ local ring,  $\rij xd$ a system
of parameters  and $J$ an ideal generated by monomials in the $x_i$, including
for each $i$ some  power of $x_i$ (e.g., if $J$ is a power of $\rij xdR$),
then $J$ is generically tightly closed.
\end{lemma}
\begin{proof} By \cite{EH}, we can write $J$  as the intersection of ideals of
the form $(x_1^{e_1},\dots,x_d^{e_d})R$, for some $e_i$. Each such ideal is
generically tightly closed by Theorem~\ref{T:norm}, whence so is $J$.
\end{proof}

\section{Local Cohomology}

Before we turn to pseudo-rationality, we must say something about local and
sheaf cohomology and their respective ultraproducts.  For
our purposes, local cohomology is most conveniently approached via
\Cech\
cohomology, which we quickly review. Let $\mathbf x:=\rij xd$ be a tuple of
elements in a Noetherian
local ring $S$ and let $\id$ be the radical of the ideal they generate. For each
$n\leq d$, define 
\begin{equation*}
\op C^n(\mathbf x;S):= \bigoplus_{1\leq i_1<i_2<\dots<i_n\leq d}
S_{x_{i_1}x_{i_2}\cdots x_{i_n}} 
\end{equation*} (with the convention that $\op C^0(\mathbf x;S)=S$). The $\op
C^n(\mathbf x;S)$ are the modules appearing in a complex $\op C^\bullet(\mathbf
x;S)$, called the \emph{algebraic \Cech\ complex} with respect to $\mathbf x$,
where the
differential $\op C^n(\mathbf x;S)\to \op C^{n+1}(\mathbf x;S)$ is given by the
inclusion maps among the localizations, with the choice of an appropriate sign
to make $\op C^\bullet(\mathbf x;S)$ a complex (see \cite[\S3.5]{BH} for more
details). The cohomology of this complex is called the \emph{local cohomology}
of $S$ with respect to $\id$ and is denoted $\op H^\bullet_\id(S)$. One shows
that $\op H^\bullet_\id(S)$ does not depend on the choice of tuple $\mathbf x$ for
which
$\op{rad}(\mathbf x S)=\id$.   If $M$ is an $S$-module, we also define the local
cohomology $\op H^\bullet_\id(M)$ of $M$  with respect to $\id$, as the
cohomology of the complex $\op C^\bullet(\mathbf x;M):=\op C^\bullet(\mathbf
x;S)\tensor_SM$. If $T$ is an $S$-algebra, then $\op H^\bullet_\id(T)=\op
H^\bullet_{\id T}(T)$. We
will be mainly interested in the top cohomology group
$\op H^d_\id(S)$ and we use the following notation. Since $\op H^d_\id(S)$ is a
homomorphic image of $\op C^d(\mathbf x;S)=S_{x_1\dots x_d}$, an arbitrary
element is the image of a fraction $\frac a{(x_1\dots x_d)^n}$
and we will denote this image by $\class a{(x_1\dots x_d)^n}S$.

For the remainder of this section, $R$ is a Noetherian local ring of equal \ch\
zero. Let $S$ be a finitely generated $R$-algebra, let $\mathbf x$ be a tuple in
$S$ and let $\id$ be the radical of $\mathbf xS$. Note that each module in the
algebraic \Cech\ complex $\op C^\bullet(\mathbf x;S)$ is  a finitely generated
$R$-algebra, whence admits an $R$-hull. The \emph{non-standard algebraic \Cech\ complex}
$\op C_\infty^\bullet(\mathbf x;S)$ over $S$ with respect to $\mathbf x$ is by
definition the complex whose $n$th module is $\rhull R{\op C^n(\mathbf x; S)}$
and for which the differentials are induced by the differentials on $\op
C^\bullet(\mathbf x;S)$.

\begin{definition}
The \emph{\ulc} of $S$ with respect
to $\id$ is by definition the cohomology of the non-standard algebraic \Cech\ complex $\op
C_\infty^\bullet(\mathbf x;S)$ and is denoted $\nslc \bullet\id S$.
\end{definition}
 
Without proof, we state that $\nslc \bullet\id S$ is independent from the
$d$-tuple generating an ideal with radical equal to $\id$. By \eqref{rhull}, the
canonical \homo{s} $\op C^n(\mathbf x; S) \to \rhull R{\op C^n(\mathbf x; S)}$
give rise to a map of complexes $\op C^\bullet(\mathbf x;S)\to\op
C_\infty^\bullet(\mathbf x;S)$, and hence for each $n\leq d$, induce  a natural
morphism 
\begin{equation*}
 j^n_\id\colon \op H^n_\id(S)\to \nslc n\id S.
\end{equation*} 
Let $\seq Sw$, $\seq \id w$ and $\seq{\mathbf x}w$ be   $R$-\sr{s} of $S$, $\id$
and $\mathbf x$ respectively. Since we can calculate the local cohomology $\op
H^\bullet_{\seq\id w}(\seq Sw)$ with aid of the algebraic \Cech\ complex of $\seq{\mathbf
x}w$ and since  taking ultraproducts commutes with cohomology, we get
\begin{equation}
\label{eq:lcul}
\nslc n\id S \iso \up w \op H^n_{\seq\id w}(\seq Sw)
\end{equation}
for each $n$. In particular, if $\varphi\colon S\to T$ is an $R$-algebra \homo\
of finite type, then the diagram
\commdiagram [lcul] {\op H^n_\id(S)} {j^n_\id} {\nslc  n\id S} {\op
H^n_\id(\varphi)}
{\nslc n\id\varphi} {\op H^n_\id(T)} {j^n_{\id T}} {\nslc n\id  T}
commutes for each $n$, where the vertical arrows are the natural maps.

\subsection*{Sheaf cohomology and local cohomology} 
Let $Y$ be a scheme,  
$\sheaf$ a quasi-coherent $\loc_Y$-module and $Z$   a closed subset of $Y$. The
collection of those global sections in $\op
H^0(Y,\sheaf)$ whose support
is contained in $Z$   is denoted $\Gamma_Z(Y,\sheaf)$ and is called the \emph{global
sections of $\sheaf$ with support in $Z$}. The derived functors of the  left-exact
functor $\Gamma_Z$ are called the \emph{cohomology with support in
$Z$} and are denoted $\op H^i_Z(Y,\sheaf)$. In case $\sheaf$ is equal to the
structure sheaf $\loc_Y$, we simply denote the cohomology with support in $Z$ by
$\op H^i_Z(Y)$.

The cohomology groups with support are connected to the usual sheaf cohomology
via an exact sequence
\begin{equation}
\label{eq:cohsupp}
\dots\to \op H^{i-1}(Y,\sheaf)\map{\rho^{i-1}} \op H^{i-1}(Y-Z,\restrict\sheaf
{Y-Z})\map{\partial^i} \op H^i_Z(Y,\sheaf)\to \op H^i(Y,\sheaf)\to\dots
\end{equation}
where $
\partial^i$ are the connecting morphisms (see for instance \cite[Corollary
1.9]{GH}). We now turn to the connection with local cohomology. 

\subsubsection*{Affine case: $Y=\op{Spec}S$}
Let $\id$ be the (radical) ideal
defining $Z$.  If we let $M:=\op H^0(Y,\sheaf)$, then $\Gamma_Z(Y,\sheaf)$ is
equal to the collection $\Gamma_\id(M)$ of all $m\in M$ which are annihilated by
some power $\id^n$. By definition, the local cohomology of $M$ with respect to 
$\id$ are the derived functors of $\Gamma_\id(M)$. Hence, for all $i$, we get an
isomorphism
\begin{equation*}
\op H^i_\id(M)\iso  \op H^i_Z(Y, \sheaf).
\end{equation*} 
  Considering
\eqref{eq:cohsupp} plus the fact
that affine schemes have no higher cohomology, we see that in the affine case 
$$
\op H^{i-1}(Y-Z,\restrict\sheaf {Y-Z})\iso \op H^i_\id(\op H^0(Y,\sheaf))
$$
 for $i>1$.

\subsubsection*{Projective case: $Y=\op{Proj}S$}
On non-affine schemes, sheaf cohomology does survive so that it has to be taken
into account. We will explain this only for 
$\sheaf:=\loc_Y$, since this is the only case we need. We fix some notation
concerning graded rings. A \emph{positively graded ring} $S$ admits a
direct sum decomposition
\begin{equation*} S=\bigoplus_{n\geq 0}\grad Sn
\end{equation*} 
where $\grad Sn$ denotes the \emph{degree $n$ part} (or \emph{$n$th homogeneous
piece}) of $S$. The irrelevant ideal of $S$ will be denoted by  $\pos S:=\oplus_{n>0}\grad
Sn$. For $x\in \pos S$, we will use $\op
D(x)$ to denote the open set of $Y=\op{Proj} S$ obtained as the complement of
the closed set $\op V(x)$, that is to say, $\op D(x):=\set{\pr\text{ homogeneous
prime of }S}{x\notin\pr}$. For a tuple $\mathbf x:=\rij xd$ with entries in $\pos
S$, we let    $\mathfrak U_{\mathbf x}$ be the
affine open covering of $Y-\op V(\mathbf xR)$ consisting of all $\op D(x_i)$.

Let  $Z$ be a closed subset of $Y$, and choose a homogeneous ideal $\id:=\mathbf
xS$ defining
$Z$ so that $\mathbf x:=\rij xd$ consists of homogeneous elements of positive
degree (for instance, let $\id:=I\pos S$ where $I$ is the homogeneous ideal
of all elements in $S$ vanishing on $Z$).  The sheaf cohomology of the
complement of $Z$, that is to say, 
$\op
H^i(Y-Z,\loc_{Y-Z})$ can be calculated as the cohomology of  the
\emph{topological
\Cech\
complex} $\op C^\bullet(\mathfrak U_{\mathbf x})$ (see \cite[III. Theorem
4.5]{Hart}).  Recall that $\op C^j(\mathfrak U_{\mathbf x})$ is the direct sum
of the 
\begin{equation*}
\op H^0(\op D(x_{i_1})\cap\dots\cap \op D(x_{i_j}),\loc_{Y-Z}) \iso \grad
{S_{x_{i_1}\dots x_{i_j}}} 0,
\end{equation*} for all $1\leq i_1<\dots<i_j\leq d$. It follows that the
topological \Cech\ complex $\op C^\bullet(\mathfrak U_{\mathbf x})$ is equal to
the degree zero part of the algebraic \Cech\ complex $\op C^\bullet(\mathbf x;S)$. Taking
cohomology (and noting that cohomology is naturally   $S$-graded) we get for
each $i=\range 2d$ that
\begin{equation}
\label{eq:ccloc}
\op H^{i-1}({Y-Z},\loc_{Y-Z}) \iso \grad {\op H^i_\id(S)} 0.
\end{equation}
This applies in particular to  $\id:=\pos S$ (whence $Z=\emptyset$) so that
\eqref{eq:ccloc} then calculates the sheaf cohomology of $Y$.
Using this isomorphism for both ${Y-Z}$ and $Y$, the natural morphism
$\rho^{i-1}\colon \op H^{i-1} (Y,\loc_Y)\to \op H^{i-1}({Y-Z},\loc_{Y-Z})$ appearing in
\eqref{eq:cohsupp} gives rise to a morphism $\grad {\op H^i_{\pos S}(S)}
0\to \grad {\op H^i_\id(S)} 0$, for each $i=\range 2d$. Let us describe this
morphism in more detail. Let $\mathbf y$ be a tuple of homogeneous elements of
positive degree so that $\pos S=\mathbf yS$. Without loss of generality, we
may assume that $\mathbf x$ is a sub-tuple of $\mathbf y$. Factoring out those
summands involving an entry not in $\mathbf x$, we get a surjective morphism of
complexes $\op C^\bullet(\mathbf y;S) \to \op C^\bullet(\mathbf x;S)$. In
cohomology, this yields morphisms $r^\bullet\colon\op H^\bullet_{\pos S}(S)\to
\op H^\bullet_\id(S)$. On the other hand, the $\rho^\bullet$ are
defined as follows. Let $\mathfrak U_{\mathbf y} -Z$ denote the collection
of
all opens $V-Z$ with $V\in\mathfrak U_{\mathbf y}$. In particular, $\mathfrak U_{\mathbf y} -Z$ is an open
covering of ${Y-Z}$. The restriction morphisms  give rise to a
morphism of topological \Cech\ complexes $\op C^\bullet(\mathfrak U_{\mathbf
y})\to \op C^\bullet(\mathfrak U_{\mathbf y}-Z)$. Since $\mathfrak
U_{\mathbf x}$ is a subcovering of $\mathfrak U_{\mathbf y}-Z$, its
topological \Cech\ complex has the same cohomology as $\op C^\bullet(\mathfrak
U_{\mathbf y}-Z)$. In other words, factoring out those summands involving an
open not in $\mathfrak U_{\mathbf x}$ yields a quasi-isomorphism of complexes
$\op C^\bullet(\mathfrak U_{\mathbf y}-Z)\to \op C^\bullet(\mathfrak 
U_{\mathbf x})$. Hence the composite map $\op C^\bullet(\mathfrak U_{\mathbf
y})\to\op C^\bullet(\mathfrak U_{\mathbf x})$ is the degree zero part of the map
$\op C^\bullet(\mathbf y;S) \to \op C^\bullet(\mathbf x;S)$. From this it is now
clear that $\rho^{i-1}$, after the identifications~\eqref{eq:ccloc}, is the 
degree zero part $\grad{r^i}0$ of $r^i$, for each $i=\range 2d$.

\subsection*{\Usc} 
Assume from now on that $S$ is a
finitely generated graded $R$-algebra, where $R$ is  an equi\ch\ zero Noetherian
local ring with trivial grading (so that in particular, each  $\grad Sn$ is an
$R$-module). Let $\seq Sw$ and $\seq{x_i}w$ be $R$-\sr{s} of $S$ and $x_i$
respectively and put $\seq{\mathbf x}w:=(\seq{x_1}w,\dots,\seq{x_d}w)$. By an
argument similar to the one in \cite[\S2.9]{SchLogTerm}, we see that almost all
$\seq Sw$ are graded and almost all $\seq{x_i}w$ are homogeneous (of degree
equal to the degree of $x_i$). Hence for each non-standard integer $\ul j:=\up
w\seq jw$ we define the \emph{degree $\ul j$ part} of $\rhull RS$ as 
\begin{equation*}
\grad {\rhull RS}{\ul j}:=\up w\grad{\seq Sw}{\seq jw}
\end{equation*}
 If we apply this to each term in the algebraic \Cech\ complex of
$\mathbf
x$ and take cohomology, we get  the  degree $\ul j$  part of the non-standard
local cohomology groups $\nslc i\id S$, and by \eqref{eq:lcul} this is also
equal to the ultraproduct of the degree $\seq jw$ parts of the local cohomology
of the \sr{s}. In view of isomorphism~\eqref{eq:ccloc}, we define for  $i=\range
2d$ the \emph{\usc} of ${Y-Z}$ as
\begin{equation*}
\nssc {i-1}{Y-Z}:= \grad {\nslc i\id S} 0.
\end{equation*} 
It follows from \eqref{eq:ccloc} that 
$$
\nssc{i-1}{Y-Z}=\up w \op H^{i-1}(\seq Yw-\seq Zw,\loc_{\seq Yw-\seq Zw}),
$$
 where
$\seq Zw:=\op V(\seq{\mathbf x}w\seq Sw)$. The natural map $j^i_\id\colon \op
H^i_\id(S)\to
\nslc i\id S$ induces in degree zero a map 
$$
u^{i-1}_{Y-Z}\colon \op
H^{i-1}({Y-Z},\loc_{Y-Z})\to \nssc {i-1}{Y-Z},
$$
fitting in a commutative diagram
\commdiagram [nssc] {\op H^{i-1}(Y,\loc_Y)} {\rho^{i-1}} {\op H^{i-1}({Y-Z},\loc_{Y-Z})}
{u^{i-1}_Y} {u^{i-1}_{Y-Z}} {\nssc {i-1}Y} {\ul\rho^{i-1}} {\nssc{i-1}{Y-Z}}
where
$\ul\rho^{i-1}$ is the ultraproduct of the restriction maps 
$$
\seq\rho w^{i-1}\colon \op
H^{i-1}(\seq Yw,\loc_{\seq Yw})\to \op H^{i-1}(\seq Yw-\seq Zw,\loc_{\seq Yw-\seq
Zw}).
$$
Indeed,
making the appropriate identifications between local cohomology and sheaf
cohomology from \eqref{eq:ccloc}, diagram~\eqref{nssc} is the degree zero part
of 
\commdiagram [rest] {\op H^i_{\pos S}(S)} {r^i} {\op H^i_\id(S)}
{j^i_{\pos S}} {j^i_\id} {\nslc i{\pos S}S} {\ul r^i} {\nslc i\id S}
where $\ul r^i$ is the  ultraproduct of the maps 
\begin{equation*}
\seq r w^i\colon \op H^i_{\pos{\seq Sw}}(\seq Sw)\to \op H^i_{\seq\id
w}(\seq Sw).
\end{equation*} It is now clear from the definition of the $\seq rw^i$ and $r^i$
that \eqref{rest} commutes, whence so does \eqref{nssc}.

\section{Pseudo-rationality}

The  notion of pseudo-rationality was introduced by  Lipman  and Teissier  to extend the notion of rational singularities to a
situation where there is not necessarily a resolution of singularities
available.

\subsection{Pseudo-rationality}\label{s:pr} 
A Noetherian local ring  $(R,\maxim)$  is called \emph{pseudo-rational}, if it
is analytically unramified, normal, \CM\ and for any projective birational map
$f\colon Y\to \op{Spec}R$ with $Y$ normal, the canonical epimorphism between the
top cohomology groups $\delta\colon H^d_\maxim(R)\to H_Z^d(Y)$ is
injective, where $Z$ is the closed fiber $\inverse f\maxim$  and $d$ the
dimension of $R$ (see \eqref{eq:cohd} below   for the   definition of
$\delta$).
Moreover, if $\op{Spec}R$ admits a desingularization
$Y\to\op{Spec}R$, then it suffices to check the above condition for just this
one $Y$ (see \cite[\S2, Remark~(a) and Example~(b)]{LT}). From
this, one can show using
 Matlis duality,  that if  $R$ is
essentially of finite type over a field  of \ch\ zero, then $R$ is
pseudo-rational \iff\ it has rational singularities. A Noetherian ring  $A$ is
called \emph{pseudo-rational}, if $A_\pr$ is pseudo-rational for every prime
ideal $\pr$ in $A$.

The key ingredient in proving Theorems~\ref{T:psrat} and
\ref{T:logterm} is the following result linking generic tight closure with
pseudo-rationality, analogous to Smith's characterization \cite{SmFrat}   in
prime \ch.

\begin{theorem}\label{T:genrat} 
If an equi\ch\ zero Noetherian local ring $R$ is
\genrat, then it is pseudo-rational.
\end{theorem}
\begin{proof} By Theorems~\ref{T:CM} and \ref{T:norm} and
Proposition~\ref{P:comp}, we know that $R$ is
analytically unramified, \CM\ and normal. Let $X:=\op{Spec} R$ and let $f\colon
Y=\op{Proj}S\to X$ be a projective birational map with $Y$ normal. In
particular, $S$ is a finitely generated graded $R$-algebra (where $R$ has the
trivial grading).  Let $i\colon R\to S$ be the
embedding identifying $R$ with $\grad S0$, let $\maxim $ be the maximal ideal of
$R$ and let $Z:=\op V(\maxim S)$ be the closed fiber of $f$.   The image of the
canonical map $\op H^d_\maxim(i)\colon \op
H^d_\maxim(R)\to \op H^d_{\maxim S}(S)$ lies entirely in degree zero whence  in
view of \eqref{eq:ccloc}, induces a morphism $\gamma^d\colon \op H^d_\maxim(R)\to
\op
H^{d-1}({Y-Z},\loc_{Y-Z})$. (Although we do
not need this, one can show that $\gamma^d$ is in fact an isomorphism, since
$\grad{\op C^\bullet(\mathbf x;S)}0=\op C^\bullet(\mathbf x;R)$.) Combining
this with   the tail of the exact sequence~\eqref{eq:cohsupp} and with
\eqref{nssc} gives a commutative diagram
\begin{equation}
\begin{aligned}
\label{eq:cohd}
\xymatrix{ 
&\op H^d_\maxim(R)\ar[d]_{\gamma^d}\ar[dr]^{\delta}\\ 
\op H^{d-1}(Y,\loc_Y) \ar[r]_-{\rho^{d-1}} \ar[d]_{u^{d-1}_Y} &\op
H^{d-1}({Y-Z},\loc_{Y-Z})\ar[d]_{u^{d-1}_{Y-Z}}\ar[r]_-{\partial^d}   & \op
H^d_Z(Y) \\
\nssc {d-1}Y \ar[r]_-{\ul\rho^{d-1}}  &\nssc{d-1}{Y-Z}
}
\end{aligned}
\end{equation}
in which the middle row is exact.

Let   $\mathbf x$ be a system of parameters in $R$ such that $\mathbf x R$ is
generically tightly closed. We need to show that the kernel of $\delta$ is zero,
hence suppose the contrary.
In particular, it must contain a non-zero element of the form $\class
ayR$, with $a\in R$ and where $y$ is the product of the entries in
$\mathbf x$.
 From  the exactness of  \eqref{eq:cohd}, we see that 
$\delta(\class ayR)=0$ means that $\gamma^d(\class ayR)$ lies in the image of
$\rho^{d-1}$.  Under the isomorphism $\op H^{d-1}({Y-Z},\loc_{Y-Z})\iso \grad{\op
H^d_{\maxim
S}(S)}0$ from \eqref{eq:ccloc}, we may identify   $\gamma^d(\class ayR)$ with
$\class ayS$. Since the square in \eqref{eq:cohd} commutes, $u^{d-1}_{Y-Z}(\class
ayS)$ lies in the image of $\ul\rho^{d-1}$. 

Let $(\seq Rw,\seq\maxim w)$   be an \sr\ of $(R,\maxim)$. By Corollary~\ref{C:srCM}, almost all $\seq Rw$ are \CM\ and
normal, whence in particular domains. Let $\seq Sw$ be an $R$-\sr\ of $S$,
put $\seq Xw:=\op{Spec}(\seq Rw)$ and $\seq Yw:=\op{Proj}(\seq Sw)$, and let
$\seq Zw:=\op V(\seq\maxim w\seq Sw)$ be the closed fiber of $\seq Yw\to\seq Xw$. Let  $\seq
aw$ and $\seq{\mathbf x}w$ be   \sr{s} of  $a$ and $\mathbf x$ respectively, and
put $\seq yw$ equal to the product of all the entries in $\seq{\mathbf x}w$. By
definition, $u^{d-1}_{Y-Z}(\class ayS)$ is the ultraproduct of the $\class{\seq
aw}{\seq yw}{\seq Sw}$. Hence by \los, almost all $\class{\seq aw}{\seq yw}{\seq
Sw}$ lie in the image of 
$$
\seq\rho w^{d-1}\colon \op H^{d-1}(\seq Yw,\loc_{\seq
Yw})\to \op H^{d-1}(\seq Yw-\seq Zw,\loc_{\seq Yw-\seq Zw})
$$  
since $\ul\rho^{d-1}$ is  the ultraproduct of the $\seq\rho w^{d-1}$.  By
the same argument as above, we have for each $w$, an exact diagram
\begin{equation}\begin{aligned}
\label{eq:cohdw}
\xymatrix{ 
&\op H^d_{\seq\maxim w}(\seq Rw)\ar[d]_{\seq\gamma w^d}\ar[dr]^{\seq \delta w}\\ 
\op H^{d-1}(\seq Yw,\loc_{\seq Yw}) \ar[r]_-{\seq\rho w^{d-1}}  &\op
H^{d-1}({\seq Yw-\seq Zw},\loc_{\seq Yw-\seq Zw}) \ar[r]_-{\seq\partial w^d}
&\op H^d_{\seq Zw}(\seq Yw).
}
\end{aligned}
\end{equation}
 By reversing the above arguments, this diagram then shows that almost each
$\class {\seq aw}{\seq yw}{\seq Rw}$ lies in
the kernel $\seq Lw$ of $\seq\delta w$. Let us briefly recall the argument from
\cite{SmFrat} how for a fixed $w$ this implies that $\seq aw$ lies in the tight
closure
of $\seq{\mathbf x}w\seq Rw$. Namely, since the Frobenius $\frob w$ acts on the
local cohomology groups,  the kernel $\seq Lw$ is invariant under its action by
functoriality. Hence 
\begin{equation}
\label{eq:frobcl}
\frob w^m (\class {\seq aw}{\seq yw}{\seq Rw}) = \class{\frob w^m(\seq aw)}
{\frob w^m(\seq yw)} {\seq Rw}\in\seq Lw.
\end{equation}
Since $\seq Lw$ is a proper subgroup of $\op H^d_{\seq\maxim w}(\seq Rw)$
(note that $\seq\delta w^d$ is non-zero),  the Matlis dual of $\seq Lw$ is   a proper
homomorphic image of the canonical module $\omega_{\seq Rw}$. Since the
canonical module  has rank one, the Matlis dual of $\seq Lw$ has   torsion,
whence so does $\seq Lw$ itself. Hence for some  non-zero 
  $\seq cw\in \seq Rw$ we have $\seq cw\seq Lw=0$. Together with
\eqref{eq:frobcl}, this yields 
\begin{equation*}
\class{\seq cw\frob w^m(\seq aw)} {\frob w^m(\seq yw)} {\seq Rw}=0
\end{equation*} 
for each $m$. Since almost each $\seq Rw$ is \CM,  we get  $\seq cw\frob
w^m(\seq aw)\in\frob w^m(\seq {\mathbf x}w)\seq Rw$, for all $m$, proving our
claim that $\seq aw$ lies in the tight closure of $\seq{\mathbf x}w\seq Rw$.
Since this holds for almost all $w$, we conclude that $a$ lies in the generic
tight closure of $\mathbf xR$, which,
 by assumption,   is just $\mathbf xR$. However, this means that $\class ayR$ is
zero, contradiction.
\end{proof}

\subsection*{Proof of Theorem~\ref{T:psrat}} 
Since all properties localize,
we may assume that $A$ and $B$ are moreover local and that $A\to B$ is a local
\homo.    Since $B$ is \genreg\ by Theorem~\ref{T:reg}, so is
$A$, by Proposition~\ref{P:cp}. Therefore, $A$ is pseudo-rational by
Theorem~\ref{T:genrat}. 
\qed

\section{Ultra-F-regular rings and log-terminal singularities}

In this section, we extend the argument from \cite{SchLogTerm} in order to prove
Theorem~\ref{T:logterm}. 

\subsection{$\mathbb Q$-Gorenstein Singularities}\label{D:Qgor} Let $R$ be an
equi\ch\ zero Noetherian local domain and put $X:=\op{Spec}R$.  We say that  
$R$ is \emph{$\mathbb Q$-Gorenstein} if it is normal and some positive multiple
of the canonical divisor $K_X$ is Cartier; the least such positive multiple is
called the \emph{index} of $R$. If $R$ is the homomorphic image of an excellent
regular local ring (which is for instance the case if $R$ is complete), then $X$
admits an \emph{embedded resolution of singularities} $f\colon Y\to X$
by \cite{Hi64}. If $E_i$ are the irreducible components of the exceptional locus
of $f$, then the canonical divisor $K_Y$ is numerically equivalent 
to $f^*(K_X)+\sum a_iE_i$ (as $\mathbb Q$-divisors), for some $a_i\in\mathbb Q$
($a_i$ is called the \emph{discrepancy} of $X$ along $E_i$; see \cite[Definition
2.22]{KM}). If all $a_i>-1$,  we call $R$ \emph{log-terminal} (in case we only
have a weak inequality, we call $R$ \emph{log-canonical}).

\subsection{Canonical cover}\label{s:kaw} Recall the construction of the
canonical cover of a $\mathbb Q$-Gorenstein local ring $R$ due to Kawamata. If 
$r$ is the index of $R$, then $\loc_X(rK_X)\iso \loc_X$, where $X:=\op{Spec}R$
and $K_X$ the canonical divisor of $X$. This isomorphism induces an $R$-algebra
structure on
\begin{equation*}
\tilde R:= H^0(X,\loc_X \oplus \loc_X(K_X)\oplus\dots\oplus \loc_X((r-1)K_X)),
\end{equation*} which is called the \emph{canonical cover} of $R$; see
\cite{Kaw}.  An important property for our purposes is  that $R\to \tilde R$ is
\'etale in codimension one (see for instance \cite[4.12]{SmVan}). We also use
the following result proven by Kawamata in \cite[Proposition 1.7]{Kaw}:

\begin{theorem}\label{T:Kaw}
Let $R$ be a homomorphic image of an equi\ch\ zero, excellent
regular local ring. If $R$ is $\mathbb Q$-Gorenstein,
then it has log-terminal singularities \iff\ its canonical cover is rational.
\end{theorem}

\begin{definition}
Inspired by Kawamata's result, we can now give a
resolution-free variant of log-terminal singularities. We call a Noetherian local ring
\emph{pseudo-log-terminal} if it is $\mathbb Q$-Gorenstein and its canonical
cover is pseudo-rational. 
\end{definition}

In the remainder of this section, $R$ is an equi\ch\
zero Noetherian local ring and  $\seq Rw$ is an \sr\ of $R$.

\subsection{Ultra-F-regularity} 
We say that  $R$ is \emph{ultra-F-regular}, if it is a domain and for each
non-zero $c\in R$, we can find an ultra-Frobenius $\ulfrob^\varepsilon$ such
that the $R$-module morphism 
\begin{equation}
\label{eq:ufc} R\to \ehull R\varepsilon\colon x\mapsto c\ulfrob^\varepsilon(x)
\end{equation}
is pure. Note  that in order for \eqref{eq:ufc} to be $R$-linear, we need to view
$\hull
R$ as an $R$-algebra via $\ulfrob^\varepsilon$, that is to say,  the target must
be taken to be $\ehull R\varepsilon$ (see \S\ref{s:uf}). Since $\hull
R=\hull{\complet R}$,   an analytically unramified local ring $R$ is
ultra-F-regular \iff\ its completion $\complet R$ is.

Over normal domains, purity and cyclical purity are the same by \cite[Theorem
2.6]{HoPure}. Hence for $R$ normal, the purity of \eqref{eq:ufc} is equivalent
to the weaker condition that  for every $x\in R$ and every ideal $I\sub  R$,
we have 
\begin{equation}\label{eq:cycpu} 
c\ulfrob^\varepsilon(x)\in\ulfrob^\varepsilon(I)\hull R
\quad\text{implies}\quad x\in I.
\end{equation}
One can show that if $R$ is moreover analytically unramified, then either
condition entails normality, and hence in that case, they are equivalent (this
follows for instance from the discussion below and the \BS\ property of generic
tight closure).

\begin{proposition}\label{P:regqreg}
 If $R$ is regular, then it is
ultra-F-regular.
\end{proposition}
\begin{proof} By the above discussion, we need only verify the weaker condition
\eqref{eq:cycpu}. In fact, we will show that for any $c$, we may take
$\varepsilon=1$ in \eqref{eq:cycpu}. Indeed, assume $c\ulfrob(x)\in\ulfrob
(I)\hull R$.
 Since $\ulfrob$ preserves regular sequences, $\ehull R1$ is a balanced big \CM\
$R$-algebra whence flat by \cite[Theorem IV.1]{SchFPD} or \cite[Lemma
2.1(d)]{HHbigCM2}. Hence 
\begin{equation*} c\in (\ulfrob(I)\hull R:\ulfrob(x))= \ulfrob(I:x)\hull R.
\end{equation*} Suppose $x\notin I$. Since $(I:x)$ then lies in the maximal
ideal of $R$, its image under $\ulfrob$ lies in the ideal of infinitesimals of
$\hull R$. Hence $\ulfrob(I:x)\hull R\cap R=(0)$,  contradicting that $c\neq 0$.
\end{proof}

\begin{theorem}\label{T:qgenreg} 
If $R$ is analytically unramified and
ultra-F-regular, then it is \genreg, whence in particular pseudo-rational.
\end{theorem}
\begin{proof} 
The last assertion follows from the first by
Theorem~\ref{T:genrat}. Since all properties are invariant under completion, we
may assume that $R$ is complete. Let
$I$ be an ideal in $R$ and $x\in\gentc I$. We want
to show that $x\in I$. By \cite[Proposition 6.24]{SchAsc}, there exists $c\in
R$ such that almost all $\seq cw$ are test elements in $\seq Rw$, where
$\seq
cw$ and $\seq Rw$ are \sr{s} of $c$ and $R$ respectively.
Let $\seq xw$
and $\seq Iw$ be \sr{s} of $x$ and $I$ respectively, so that almost all $\seq
xw\in\seq Iw^*$. Hence, for almost all $w$ and all $e$, we have
\begin{equation}
\label{eq:cwe}
\seq cw\frob w^e(\seq xw)\in\frob w^e(\seq Iw)\seq Rw.
\end{equation}
By assumption, there is an ultra-Frobenius $\ulfrob^\varepsilon$ so that
$x\mapsto c\ulfrob^\varepsilon(x)$ is pure whence cyclically pure, that is to
say, so that \eqref{eq:cycpu} holds. Let $\varepsilon$ be the ultraproduct of
integers $\seq ew$. Taking $e$ equal to $\seq ew$ in \eqref{eq:cwe} and taking
ultraproducts shows that $c\ulfrob^\varepsilon(x)\in \ulfrob^\varepsilon(I)\hull
R$. Therefore,  from \eqref{eq:cycpu} we get $x\in I$, as we
wanted to show.
\end{proof}

\begin{proposition}\label{P:et} 
Let $R\sub  S$ be a finite extension of
Noetherian local domains which is \'etale in codimension one. Let $c$ be a
non-zero element of $R$ and $\ulfrob^\varepsilon$ an ultra-Frobenius. If $R\to
\ehull R\varepsilon\colon x\mapsto c\ulfrob^\varepsilon(x)$ is pure, then so is 
its base change $S\to \ehull S\varepsilon\colon x\mapsto
c\ulfrob^\varepsilon(x)$. 

In particular, if $R$ is ultra-F-regular, then so is $S$.
\end{proposition}
\begin{proof}
 Let $R\sub  S$ be an arbitrary finite extension  of  $d$-dimensional
Noetherian local  domains and fix a non-zero element $c\in R$ and an
ultra-Frobenius  $\ulfrob^\varepsilon$. Let $\mathfrak n$ be the maximal ideal
of $S$ and $\omega_S$ its canonical module. I claim that if $R\sub  S$ is
\'etale, then 
\begin{equation}
\label{eq:etul} \ehull S\varepsilon\iso S\tensor_R\ehull R\varepsilon.
\end{equation}
Assuming the claim, let $R\sub  S$ now only be \'etale in codimension one.  It
follows from the claim that the supports of the kernel and the cokernel of the
natural map $S\tensor_R\ehull R\varepsilon\to \ehull S\varepsilon$ have 
codimension at least two. Hence the same is true for the base change
\begin{equation*}
\omega_S\tensor_R\ehull R\varepsilon\to \omega_S\tensor_S\ehull S\varepsilon.
\end{equation*} Applying the top local cohomology functor $\op H_{\mathfrak
n}^d$, we get from the long exact sequence of local cohomology and  Grothendieck
Vanishing, an isomorphism
\begin{equation}
\label{eq:top}
\op H_{\mathfrak n}^d(\omega_S\tensor_R\ehull R\varepsilon) \iso \op
H_{\mathfrak n}^d(\omega_S\tensor_S\ehull S\varepsilon).
\end{equation}
Recall that by  Grothendieck duality, $\op H_{\mathfrak n}^d(\omega_S)$ is the
injective hull $E$ of the  residue field of $S$. 

Let $c_{\varepsilon,R}$ denote the $R$-linear morphisms $R\to \ehull
R\varepsilon\colon x\mapsto c\ulfrob^\varepsilon(x)$. For an arbitrary
$R$-module $M$, let  $c_{\varepsilon,R,M}\colon M\to M\tensor_R \ehull
R\varepsilon$ be the base change of $c_{\varepsilon,R}$ over $M$. In particular,
we have a commutative diagram
\begin{equation*}
\CD
 \omega_S @>c_{\varepsilon,R,\omega_S}>>  \omega_S\tensor_R\ehull
R\varepsilon\\
@| @VVV  \\
\omega_S @>>c_{\varepsilon,S,\omega_S}>  \omega_S\tensor_S\ehull S\varepsilon.
\endCD
\end{equation*}
Taking top local cohomology yields  the
outer square in the following   commutative diagram

\begin{equation}
\label{E}
\CD E=\op H_{\mathfrak n}^d(\omega_S) @>c_{\varepsilon,R,E}>> E\tensor_R \ehull
R\varepsilon@>>> \op H_{\mathfrak n}^d(\omega_S\tensor_R\ehull R\varepsilon) \\
@| @VVV @VV\iso V \\
E=\op H_{\mathfrak n}^d(\omega_S) @>>c_{\varepsilon,S,E}> E\tensor_S\ehull
S\varepsilon @>>> \op H_{\mathfrak n}^d(\omega_S\tensor_S\ehull S\varepsilon)
\endCD
\end{equation}
where the isomorphism at the right comes from \eqref{eq:top}. Since
$c_{\varepsilon,R}$ is pure, so is its base change $c_{\varepsilon,R,\omega_S}$.
Purity is preserved when taking cohomology, so that the top composite map in
\eqref{E} is pure, whence so is the bottom composite map, since it is isomorphic
to it. Since $c_{\varepsilon,S,E}$ is a factor of this map, it is itself pure,
whence in particular injective.  By \cite[Lemma 2.1(e)]{HHbigCM2}, to verify the
purity of $c_{\varepsilon, S}$, one only needs to show that its base change
$c_{\varepsilon,S,E}$ over $E$ is injective, and this is exactly what we just
showed.

To prove the claim \eqref{eq:etul}, observe that if $R\to S$ is \'etale with \sr\
$\seq Rw\to\seq Sw$, then almost all  of these are \'etale. Indeed, by  
\cite[Corollary 3.16]{Milne}, we can write $S$ as $\pol RX/I$, with $X=\rij Xn$
and $I=\rij fn\pol RX$, such that the Jacobian $\op J(f_1,\dots,f_n)$   is a
unit in $R$, and by \los, this property is preserved for almost all \sr{s}.
Quite generally,
if $C \to D$ is an \'etale extension of \ch\ $p$ domains, then we have
for each $e$ an isomorphism $\mathstrut^eD\iso D\tensor_C\mathstrut^eC$ (see
for
instance \cite[p. 50]{HuTC} or the proof of \cite[Theorem 4.15]{SmVan}). Applied
to the current situation, we get $\mathstrut^e\seq Sw\iso\seq Sw\tensor_{\seq Rw}\mathstrut^e\seq
Rw$  (see \cite[p. 50]{HuTC}). Therefore,  applied with
$e=:\seq ew$, where $\seq ew$ is an \sr\ of $\varepsilon$, we get after taking
ultraproducts,  
$$
\ehull S\varepsilon\iso \hull S\tensor_{\hull R} \ehull R\varepsilon\iso S\tensor_R\ehull R\varepsilon
$$
as required, where we used the isomorphism $\hull S\iso S\tensor_R\hull R$, which
holds by \cite[\S4.10.4]{SchAsc}, since $R\to S$ is finite.

To prove the last assertion, we have to show that we can find for each non-zero
$c\in S$ an ultra-Frobenius $\ulfrob^\varepsilon$ such that $c_{\varepsilon,S}$ 
is pure. However, if we can do this for some non-zero multiple of $c$, then we
can also do this for $c$, and hence, since $S$ is finite over $R$,  we may
assume without loss of generality that $c\in R$. Since $R$ is ultra-F-regular,
we can find therefore an ultra-Frobenius $\ulfrob^\varepsilon$ such that
$c_{\varepsilon,R}$  is pure, and hence by the first assertion, so is then
$c_{\varepsilon,S}$, proving that $S$ is ultra-F-regular. 
\end{proof}

\begin{proposition}\label{P:cpqreg}
Let $R\to S$ be a cyclically pure \homo\ of equi\ch\ zero Noetherian local
rings. If $S$ is ultra-F-regular and analytically unramified, then so is $R$.
\end{proposition}
\begin{proof}
Since $\complet R\to \complet S$ is again cyclically pure by \cite[Lemma
6.7]{SchAsc}, we may assume without
loss of generality  that $S$ is complete. Let $c\in R$ be non-zero 
and let $\ulfrob^\varepsilon$ be an ultra-Frobenius for which the $S$-module
morphism
\begin{equation}
\label{eq:spur}
  c_{\varepsilon,S}\colon S\to \ehull S\varepsilon\colon x\mapsto
c\ulfrob^\varepsilon(x)
\end{equation}
is pure. We want to show that the same is true upon replacing $S$ by $R$, that
is to say, that $c_{\varepsilon,R}$  is pure.
Since $S$ is
\genreg\
by Theorem~\ref{T:qgenreg},
 so is $R$ by Proposition~\ref{P:cp}. Hence $R$ is in particular normal by
Theorem~\ref{T:norm}, so that it suffices to  verify \eqref{eq:cycpu}. Let
$x\in R$ and $I\sub R$ be such that
$c\ulfrob^\varepsilon(x)\in\ulfrob^\varepsilon(I)\hull R$. Therefore, $x$ 
belongs to $IS$ by \eqref{eq:spur},
 whence to  $IS\cap R=I$ by cyclical purity.
\end{proof}

Note that in the proof,   the condition that $S$ is analytically
unramified was only used to get the normality of $R$. 

\subsection*{Proof of Theorem~\ref{T:logterm}} 
Proposition~\ref{P:regqreg} yields that $S$ is ultra-F-regular, whence so
is $R$ by Proposition~\ref{P:cpqreg}. Let $\tilde R$ be the canonical cover 
of $R$. By Proposition~\ref{P:et}, also $\tilde R$ is ultra-F-regular, whence
pseudo-rational by Theorem~\ref{T:qgenreg}. 
\qed


\providecommand{\bysame}{\leavevmode\hbox to3em{\hrulefill}\thinspace}
\providecommand{\MR}{\relax\ifhmode\unskip\space\fi MR }
\providecommand{\MRhref}[2]{%
  \href{http://www.ams.org/mathscinet-getitem?mr=#1}{#2}
}
\providecommand{\href}[2]{#2}

\end{document}